\documentclass[12pt,letterpaper]{amsart}

\newcounter{myenum}

\newcommand{\indentalign}{\hspace{0.3in}&\hspace{-0.3in}}
\newcommand{\la}{\langle}
\newcommand{\ra}{\rangle}

\renewcommand{\Im}{\operatorname{Im}}

\usepackage{amssymb}
\usepackage{amsthm}
\usepackage{amsxtra}
\usepackage{graphicx}

\newtheorem{theorem}{Theorem}

\newtheorem{lemma}[theorem]{Lemma}
\newtheorem{corollary}[theorem]{Corollary}

\theoremstyle{remark}
\newtheorem{remark}[theorem]{Remark}

\numberwithin{equation}{section}
\numberwithin{theorem}{section}

\ifx\pdfoutput\undefined
  \DeclareGraphicsExtensions{.pstex, .eps}
\else
  \ifx\pdfoutput\relax
    \DeclareGraphicsExtensions{.pstex, .eps}
  \else
    \ifnum\pdfoutput>0
      \DeclareGraphicsExtensions{.pdf}
    \else
      \DeclareGraphicsExtensions{.pstex, .eps}
    \fi
  \fi
\fi

\title[Asymptotically linear solutions to NLS and Hartree]
{Asymptotically linear solutions in $H^1$ of the 2-d defocusing nonlinear Schr\"odinger and Hartree equations}

\author{Justin Holmer}
\address{University of California, Berkeley}
\author{Nikolaos Tzirakis}
\address{University of Illinois, Urbana-Champaign}
\begin{document}

\maketitle

\begin{abstract}
In the 2-d setting, given an $H^1$ solution $v(t)$ to the linear Schr\"odinger equation $i\partial_t v +\Delta v =0$, we prove the existence (but not uniqueness) of an $H^1$ solution $u(t)$ to the defocusing nonlinear Schr\"odinger (NLS) equation $i\partial_t u + \Delta u -|u|^{p-1}u=0$ for nonlinear powers $2<p<3$ and the existence of an $H^1$ solution $u(t)$ to the defocusing Hartree equation $i\partial_t u + \Delta u -(|x|^{-\gamma}\star|u|^{2})u=0$ for interaction powers $1<\gamma<2$, such that $\|u(t)-v(t)\|_{H^1} \to 0$ as $t\to +\infty$.  This is a partial result toward the existence of well-defined continuous wave operators $H^1 \to H^1$ for these equations.  For NLS in 2-d, such wave operators are known to exist for $p\geq 3$, while for $p\leq 2$ it is known that they cannot exist.  The Hartree equation in 2-d only makes sense for $0<\gamma<2$, and it was previously known that wave operators cannot exist for $0<\gamma\leq 1$, while no result was previously known in the range $1<\gamma<2$.  Our proof in the case of NLS applies a new estimate of Colliander-Grillakis-Tzirakis \cite{CGT} to a strategy devised by Nakanishi \cite{N}. For the Hartree equation, we prove a new correlation estimate following the method of \cite{CGT}. 
\end{abstract}

\section{Introduction}
In this paper we consider the asymptotic behavior of $H^1$ solutions on $\mathbb{R}^n$ to the following two equations:  the defocusing nonlinear Schr\"odinger equation (NLS)
\begin{equation}\label{nls}
i\partial_tu + \Delta u -|u|^{p-1}u=0, 
\end{equation}
with nonlinear exponents $1<p<(n+2)/(n-2)$, and the defocusing Hartree equation
\begin{equation}\label{hnls}
i\partial_t u+ \Delta u -(|x|^{-\gamma} \star |u|^{2})u=0, 
\end{equation}
with interaction exponents $0<\gamma <\min(n,4)$.  For both equations, there is a local well-posedness theory (see  \cite{tc}\cite{tt} for exposition) for the initial-value problem yielding solutions $u\in C(t\in(-T,T); H_x^1(\mathbb{R}^n))$, $T=T(\|u_0\|_{H^1})$ that satisfy mass conservation:
$$
\|u(t)\|_{L^{2}}=\|u_{0}\|_{L^{2}}
$$
and energy conservation:
$$E(u)(t)=\frac{1}{2}\int |\nabla u(t)|^{2}dx+\frac{1}{p+1}\int |u(t)|^{p+1}dx=E(u_{0}) $$
for NLS (1.1) and
$$
E(u)(t)=\frac{1}{2}\int |\nabla u(t)|^{2}dx+\frac{1}{4}\int (|x|^{-\gamma} \star |u(t)|^{2})|u(t)|^2dx=E(u_{0}). 
$$
for Hartree \eqref{hnls}.

The conservation laws provide an $H^1$ \emph{a priori} bound on the solutions guaranteeing that the local-in-time solutions are in fact global in $H^1$.  In this situation, we seek to describe the asymptotic behavior of solutions in time.   The comparison between the dynamics of the nonlinear equations and their linear counterparts gives rise to the following two questions.  Fix a data class $X$ (defined by some spatial norm).  We note that since $u(x,t)$ solves NLS/Hartree $\implies$ $\bar u(x,-t)$ solves NLS/Hartree, studying asymptotics at $+\infty$ is equivalent to studying them at $-\infty$.  (But of course, a given nonlinear solution $u(t)$ will in general have a different limiting behavior at $+\infty$ than at $-\infty$.)  Let $U(t)=e^{it\Delta}$ denote the linear Schr\"odinger evolution group.

\medskip

\noindent \emph{Existence of wave operators}.  Given a linear solution $v_+(t)=U(t)u_+$ with $u_+\in X$, does there exists a global nonlinear solution $u(t)$ to NLS or Hartree such that $U(-t)u(t)\in X$ for all $t$ and 
\begin{equation}
\label{E:convergence}
\|U(-t)u(t)-u_+ \|_X \to 0 
\end{equation}
as $t\to +\infty$?\footnote{For some spaces $X$, such as $X=\Sigma$, the linear flow is not unitary and it is not known whether the statement $\|U(-t)u(t)-u_+\|_X\to 0$ is equivalent to $\|u(t)-U(t)u_+\|_X \to 0$.  On this question see Begout \cite{pb}. In this paper, we consider $X=H^1$ on which the linear flow is unitary.}   Is the solution $u(t)$ the unique solution satisfying \eqref{E:convergence} in a certain space-time function class $Y$?  If so, then we can define an operator $\Omega_+:X\to X$ sending $u_+\mapsto u(0)$, called the \emph{wave operator}.  Is $\Omega_+$ continuous? 

\medskip

\noindent \emph{Asymptotic completeness}. Given a global solution $u(t)$ to NLS or Hartree such that $U(-t)u(t)\in X$ for all $t\in \mathbb{R}$, does there exist $u_+\in X$ such that \eqref{E:convergence} holds?  Such a $u_+$ is necessarily unique.
 
\medskip

In this paper, we consider the above problems for $X=H^{1}$.  There is a large amount of literature dealing with the weighted space $\Sigma$, where $\|\phi\|_\Sigma = \|xu\|_{L^2}+ \|u\|_{H^1}$.  However, to remain focused in our discussion, we will not mention these results here.  Also, some results are available for the focusing analogues ($+$ sign in front of the nonlinearities), but we will also not discuss this case either.

The analysis of these problems breaks into different ranges of $p$ and $\gamma$.  The NLS equation has scaling $u(x,t) \mapsto \lambda^\frac{2}{p-1}u(\lambda x,\lambda^2 t)$ and the Hartree equation has scaling $u(x,t) \mapsto \lambda^\frac{n-\gamma+2}{2}u(\lambda x,\lambda^2t)$.  The $H^1$ subcritical restrictions are thus $1<p<(n+2)/(n-2)$ and $0<\gamma<\min(n,4)$.  The $L^2$ critical exponent for NLS is $p=1+\frac4n$ and for Hartree is $\gamma=2$.  The basic heuristic is that for $L^2$ subcritical exponents one has good short-time control, for $L^2$ supercritical exponents one has good long-time control, and in the $L^2$ critical case the scaling shows that short times and long times are equivalent.  The expectation therefore is that the questions of asymptotic completeness and existence of wave operators should be easier in the $L^2$ supercritical case (although not necessarily false in the $L^2$ subcritical case).  It turns out that an important boundary is $p_0=1+\frac2n$ for NLS and $\gamma_0=1$ for Hartree.  For nonzero solutions $u(t)$ to NLS with $p\leq p_0$ or Hartree with $\gamma\leq \gamma_0$,  $U(-t)u(t)$ does not converge in $L^2$ (proved by Strauss \cite{ws}, Barab \cite{Barab}, and Tsutsumi \cite{Tsut85} for NLS and by Hayashi-Tsutsumi \cite{HT} and Glassey \cite{G77} for Hartree).  Thus wave operators cannot exist in any reasonable sense, and in this regime the positive results involve the construction of \emph{modified} wave operators, a topic we do not discuss here.  The region $p<p_0$ and $\gamma<\gamma_0$ is called long-range and the region $p>p_0$ and $\gamma>\gamma_0$ is called short-range.

In the $L^2$ supercritical case $p>1+\frac4{n}$, $\gamma>2$, a fixed point argument using the Strichartz estimates (similar to the one used to solve the local Cauchy problem) establishes both the existence of continuous wave operators and asymptotic completeness for small initial data.
To extend the asymptotic completeness results to large initial data, one usually proceeds through the derivation of \emph{a priori} estimates for general solutions to furnish the needed decay estimates as $t\to \infty$.  These estimates take advantage of the momentum conservation law
$$
\vec{p}(t)=\Im \int_{\Bbb R^n}\bar{u}\nabla u dx.
$$
We can establish, for example, the generalized virial identity (Lin-Strauss \cite{ls}), valid for $n \geq 3$,
\begin{equation}\label{linstr}
\int_{0}^{T}\int_{\Bbb R^n}(-\Delta \Delta |x|)|u(x,t)|^2dxdt+\frac{2(p-1)}{p+1}\int_{0}^{T}\int_{\Bbb R^n}\frac{|u(x,t)|^{p+1}}{|x|}dxdt \lesssim
 \sup_{[0,T]}|M_{a}(t)|
\end{equation}
where $u$ is a solution to (1.1) and $M_{a}(t)$ is the Morawetz action defined by
\begin{equation}\label{mora}
M_{a}(t)=2\int_{\Bbb R^n}\nabla a(x) \cdot \Im(\bar{u}(x)\nabla u(x))dx.
\end{equation}
Using this point of view, one can establish \emph{a priori} estimates that completely answer the questions above in the $L^2$ supercritical case $p>1+\frac{4}{n}$.  In the case of NLS, Ginibre-Velo \cite{gv1} proved asymptotic completeness for $n \geq 3$ using \eqref{linstr} and subsequently, Nakanishi \cite{kn} settled the case of low dimensions $n=1,2$ using special modifications of \eqref{linstr}.  These results are very technical but nowadays we can simplify the proofs using correlation estimates -- see Colliander-Grillakis-Tzirakis \cite{CGT} and the references therein.  Corresponding results for Hartree were obtained in Ginibre-Velo \cite{GV} and Nakanishi \cite{Nak}.

As we have mentioned, in the $L^2$ supercritical case, the existence of continuous wave operators can be established by a fixed point argument using the Strichartz estimates.   However, for a given $u_+$, the corresponding $u(t)$ is only uniquely given (by this argument) in $Y=L_t^qL_x^r$, a Strichartz space  -- this is a statement of \emph{conditional uniqueness}.  However, the above global decay estimates plus the use of the Strichartz estimates can be used to establish that any solution to NLS or Hartree in $C(\mathbb{R}; H_x^1)$ must also belong to $Y$ thereby establishing \emph{unconditional uniqueness}.   We point out, however, that conditional uniqueness is adequate to give well-defined wave operators and conditional as opposed to unconditional uniqueness statements are anyway quite standard in low regularity well-posedness theory.

In the $L^2$ critical case, one can again establish the existence of continuous wave operators with conditional uniqueness by a direct fixed point argument using the Strichartz estimates.  However, the Morawetz global decay estimates are alone insufficient to establish unconditional uniqueness for $X=H^1$.   Asymptotic completeness for $X=H^1$ is an interesting and difficult problem.  In the case of NLS, Killip-Tao-Visan-Zhang \cite{tv}\cite{ktv} recently showed asymptotic completeness for $n \geq 2$ assuming radial initial data (their result is true even for $L^2$ radial data).  Following their work, Miao-Xu-Zhao \cite{MXZ} established the analogous result for radial solutions to Hartree with $\gamma=2$, $n\geq 3$.

The $L^2$ subcritical case $p_0=1+\frac{2}{n}<p<1+\frac{4}{n}$ and $\gamma_0=1<\gamma<2$ in $X=H^1$ is a rather open area.  There are satisfactory results available for $X=\Sigma$ that we shall not discuss here; see Cazenave \cite{tc} Chapter 7 for an overview and references.\footnote{In weighted spaces, the existence of continuous wave operators is known for the whole eligible $L^2$ subcritical range and asymptotic completeness is known to hold for $p>p_S$ and $\gamma>\gamma_S$, where $p_S>p_0$ and $\gamma_S>\gamma_0$ are certain exponents that emerge from the use of the pseudoconformal conservation law.  The main remaining open problem (as far as we know) in this setting is to determine whether or not asymptotic completeness holds for $p_0<p<p_S$ and $\gamma_0<\gamma<\gamma_S$.}  Regarding the energy space $X=H^1$, a result that essentially goes back to Segal \cite{is} states that for NLS, one has \emph{weak} asymptotic completeness.  This means that for any solution $u(t)$ to NLS, $U(-t)u(t)$ converges weakly in $H^1$ as $t\rightarrow \pm \infty$.  Presumably this result carries over to Hartree.  The obvious drawback of this result is that the weak topology cannot even distinguish standing waves from the zero solution.   Nakanishi \cite{N} showed, by a compactness argument for dimensions $n \geq 3$ and $X=L^{2}$ or $X=H^{1}$, a result in the direction of the existence of wave operators:  Given any solution $v_+\in C(\Bbb R ; X)$ of the linear Schr\"odinger equation, there exists a (strong) solution  $u \in C(\Bbb R ; X)$ of NLS such that $\lim_{t \rightarrow \infty}\|u(t)-v_+(t)\|_{X}=0$.  The uniqueness (even in the conditional sense) of such a solution would give well-defined wave operators in $X$, but it is unknown if uniqueness holds or not.  Nakanishi's proof adapts to establish an analogous result for Hartree for $n\geq 3$.\footnote{K. Nakanishi, personal communication.}

In this paper we try to complete the picture in low dimensions and follow closely the argument in Nakanishi \cite{N}.   We have the following theorem:

\begin{theorem}
\label{Theorem1}
Let $n=2$ and $2<p<3$.  Then, for any solution $v\in C(\Bbb R ; H^1)$ of the linear Schr\"odinger equation, there exists a (strong) solution $u \in C(\Bbb R ; H^1)$ of NLS satisfying
$$\lim_{t \rightarrow \infty}\|u(t)-v(t)\|_{H^1}=0.$$
The same result holds for Hartree in the case $n=2$ and $1<\gamma<2$.
\end{theorem}

Nakanishi's proof is restricted to dimensions 3 and higher because 3 is the minimum dimension where the time decay of the linear dispersive estimate is integrable.  In two dimensions, the decay just fails to be integrable, but we can still obtain the result by making use of a gain introduced by the global \emph{a priori} Morawetz estimates.  In the case of NLS, the needed estimate appears in Colliander-Grillakis-Tzirakis \cite{CGT}
\footnote{For the case of the nonlinear Schr\"odinger equation these estimates have been obtained independently and simoultaneously by Planchon and Vega \cite{pv}.}, and for Hartree we supply here in \S \ref{S:Morawetz} a proof of an analogous estimate.   We note that, like Nakanishi, we do not know how to prove a uniqueness statement that would give the existence of a well-defined wave operator.  Presumably, a uniqueness proof would adapt to establish continuity of the wave operator.  In addition we cannot extend the result to dimension one due to the very weak decay of the linear solution, even though Morawetz estimates are now available in that setting. 

\medskip

\noindent{\bf Acknowledgement.} We thank Kenji Nakanishi for useful discussions and for pointing his work in \cite{N} to us. We also thank the anonymous referees for their comments and suggestions.

\section{Global decay estimates}
\label{S:Morawetz}

In this section, we discuss the following theorem.
 
\begin{theorem}[Correlation estimate for NLS and Hartree]
\label{theorem3}
 Let $u$ be an $H^{\frac{1}{2}}$ solution to NLS \eqref{nls} or Hartree \eqref{hnls} on the space-time slab $I\times {\Bbb R}^n$, $n \geq 2$. Then
\begin{equation}\label{22d}
\|D_x^{-\frac{n-3}{2}}(|u|^2)\|_{L_I^2L_{x}^{2}} \lesssim \|u\|_{L_I^{\infty}\dot{H}_{x}^{\frac{1}{2}}}\|u\|_{L_I^{\infty}L_x^2}
\end{equation}
\end{theorem}

The proof of Theorem \ref{theorem3} in the case of NLS can be found in \cite{CGT}.  Here, we prove the Hartree case, following the method developed in \cite{CGT}.  
We will prove in details the 2d estimate. The computations for the general case are almost identical. The derivation of the estimate with the Hartree nonlinearity is more complicated because the time derivative of the pressure is not any more a divergence of a vector field.  On the other hand, we still can view the evolution equation as describing the evolution of a compressible dispersive fluid whose pressure is a function of the density. Morally, this fact and the defocusing character of the nonlinearity is enough to guarantee the existence of positive terms in the expansion of the time derivative of the Morawetz action. Note in addition that if $u_1(x_{1},t)$ and $u_2(x_{2},t)$ are solutions of $iu_t+\Delta u=f(u)$, then the tensor product $u(x,t)=(u_1\otimes u_2)(x,t)=u_1(x_1,t)u_2(x_2,t)$ satisfies the same equation with nonlinearity 
$$f(x,t)=f(u_1(x_1,t))\otimes u_2(x_2,t)+f(u_2(x_2,t))\otimes u_1(x_1,t) \, .$$ 
Here, $x=(x_1,x_2)$ and $\Delta =\Delta_{x_1}+\Delta_{x_2}$.  Sometimes we write $x \in \Bbb R^n \otimes \Bbb R^n$ and we mean $x=(x_1,x_2)$ where $x_1 \in \Bbb R^n$ and $x_2 \in \Bbb R^n$. In particular if the equation is defocusing, the tensor product also satisfies a defocusing equation.  We can thus derive correlation estimates by applying the Morawetz and Lin-Strauss method to tensor products of solutions. 

\begin{remark}
The reader will notice that we use Einstein's summation convention throughout the paper. According to this convention, when an index variable appears twice in a single term, once in an upper (superscript) and once in a lower (subscript) position, it implies that we are summing over all of its possible values.
\end{remark}

\begin{proof}
The main technique is to use commutator vector operators acting on the local conservation laws of the equation.  We define 
$$M_{a}^{\otimes_{2}}(t)= \int_{\Bbb R^{n}\otimes \Bbb R^{n}}\nabla a(x)
\cdot \Im\left(\overline{u_1\otimes u_2}(x)\nabla (u_1\otimes u_2(x))\right )dx$$
which is the Morawetz action for the tensor product of two solutions $u:=(u_{1}\otimes u_{2})(x,t)$ where $x=(x_1,x_2) \in \Bbb R^n \otimes \Bbb R^n$. 
If we specialize to the case $u_1=u_2$, $a(x)=a(x_1,x_2)=|x_1-x_2|$, $n=2$, and observe that 
$$\partial_{x_1}a(x_1,x_2)=\frac{x_1-x_2}{|x_1-x_2|}=-\frac{x_2-x_1}{|x_1-x_2|}=-\partial_{x_2}a(x_1,x_2) \, ,$$ 
we can view $M(t) := M_{a}^{\otimes_{2}}(t)$ as
\begin{equation}
\label{commute}
M(t)=2\int_{\Bbb R^2 \otimes \Bbb R^2}\frac{x_1-x_2}{|x_1-x_2|}\cdot \left( \vec{p}(x_1,t)\rho(x_2,t)-\vec{p}(x_2,t)\rho(x_1,t)\right) \, dx_1dx_2 \,,
\end{equation}
where $\rho=\frac{1}{2}|u|^2$ is the mass density and $p_{j}=\Im(\bar{u}\partial_{j}u)$ is the momentum density.  Now define the integral operator
$$D^{-1}f(x)=\int_{\Bbb R^2}\frac{1}{|x-y|}f(y)dy \, ,$$ 
where $D$ stands for the derivative. This is indeed justified because for $n=2$ the distributional Fourier transform of $\frac{1}{|x|}$ is $\frac{1}{|\xi|}$. The main observation is that we can write the action term $M(t)$ using a commutator in the following manner:
\begin{equation}
\label{E:M}
M(t)=2\langle [x;D^{-1}]\rho(t) \  | \ \vec{p}(t)\rangle.
\end{equation}
This equation follows from an elementary rearrangement of the terms of \eqref{commute}.  This suggests that the estimate \eqref{22d} is derived using the vector operator $\vec{X}$ defined by
\begin{equation}
\label{E:X}
\vec{X}=[x;D^{-1}] \, .
\end{equation}
We change notation and write $x_1:=x$ and $x_2:=y$.  The crucial property is that the derivatives of this operator $\partial_{j}X^{k}$ form a positive definite operator.  Note that in physical space
\begin{equation}
\label{E:Xint}
\vec{X}f(x)=\int_{\Bbb R^2}\frac{x-y}{|x-y|}f(y)dy \, ,
\end{equation}
and a calculation shows that

\begin{equation}
\label{E:dX}
(\partial_{j}X^{k})f(x)=\int_{\Bbb R^2}\eta_{j}^{k}(x,y)f(y)dy \, ,
\end{equation}
where
$$\eta_{j}^{k}(x,y)=\frac{\delta_{j}^{k}}{|x-y|}-\frac{(x_j-y_j)(x^k-y^k)}{|x-y|^3} \, .$$
Thus, $\partial_j X^k$ is a positive definite operator, meaning that for any $\alpha(x)=(\alpha_1(x),\alpha_2(x))$, we have
$$\int \partial_j X^kf(x) \; \alpha^j(x) \alpha_k(x) \, dx = \int \beta(y) \,f(y) \, dy \quad \text{with} \quad \beta\geq 0 \,.$$
Note also from \eqref{E:dX} that the divergence of the vector field $\vec{X}$ is given by
$$\nabla \cdot \vec{X}=D^{-1}.$$
Now, if we differentiate \eqref{E:M} (using \eqref{E:X}), we obtain that
\begin{equation}\label{dermor}
\frac{1}{2}\partial_{t}M(t)=\langle \vec{X}\partial_{t}\rho(t) \  | \ \vec{p}(t)\rangle-\langle \vec{X}\cdot \partial_{t}\vec{p}(t) \  | \ \rho(t)\rangle \, ,
\end{equation}
where we have used the fact that $\vec{X}$ is an antisymmetric operator.  Now recall the local conservation of mass
\begin{equation}\label{localmass}
\partial_{t}\rho+\partial_{j}p^{j}=0 \,.
\end{equation}
Set $V(x)=|x|^{-\gamma}$.  An explicit calculation shows that if we differentiate in time the pressure $p_{k}$, we obtain 
\begin{equation}
\label{localmom}
\partial_{t}p_k+\partial_{j}\big( \sigma_{k}^{j}+\delta_{k}^{j}\left( -\Delta \rho \right)\big)+4\rho \; \partial_{k}( V \star \rho)=0 \, ,
\end{equation}
where
\begin{equation}
\label{E:sigma}
\sigma_{k}^{j}=\frac{1}{\rho}(p^jp_k+\partial^j \rho \partial_k \rho) \, .
\end{equation}
Applying the operator $X$ to the equation \eqref{localmass} and contracting with $p_{k}$, and similarly applying the operator $X$ to equation \eqref{localmom} and contracting with $\rho$, we obtain from \eqref{dermor} that
\begin{align*}
\frac{1}{2}\partial_{t}M &=
\begin{aligned}[t] 
& \langle \sigma_{k}^{j} \  | \ (\partial_{j}X^{k})\rho \rangle -\langle p^{j} \  | \ (\partial_{j}X^{k})p_k\rangle  \\
&- \langle \Delta \rho \  | \ (\partial_{j}X^{j})\rho\rangle+4\langle X^{k}[ \rho \, \partial_{k}(V \star \rho)] \ | \ \rho \rangle
\end{aligned} 
\\
&= 
\begin{aligned}[t]
&\langle \sigma_{k}^{j} \  | \ (\partial_{j}X^{k})\rho\rangle - \langle p^{j} \  | \ (\partial_{j}X^{k})p_k\rangle \\
&-\langle \Delta \rho \  | \ (\partial_{j}X^{j})\rho\rangle  -4\langle \rho \, \partial_{k}(V \star \rho) \ | \ X^{k}\rho \rangle \,,
\end{aligned}
\end{align*}
where we used again the fact that $\vec{X}$ is antisymmetric.
Substituting \eqref{E:sigma}, we have that
$$\frac{1}{2}\partial_{t}M=P_{1}+P_{2}+P_{3}+P_{4}$$
where
\begin{align*}
& P_1:=\langle \rho^{-1}\partial_{k}\rho \partial^{j}\rho \  | \ (\partial_{j}X^{k})\rho\rangle \\
& P_2:=\langle \rho^{-1}p_kp^j \  | \ (\partial_{j}X^{k})\rho\rangle-\langle p^j \  | \ (\partial_{j}X^{k})p_k\rangle \\
& P_3:=\langle (-\Delta \rho) \  | \ (\partial_{j}X^{j})\rho\rangle=\langle (-\Delta \rho) \  | \ (\nabla \cdot \vec{X})\rho\rangle \\
& P_4:=-4\langle \rho \, \partial_{k}(V \star \rho) \ | \ X^{k}\rho \rangle 
\end{align*}
The term $P_1$ is clearly positive since $\partial_{j}X^{k}$ is a positive definite operator. Let's analyze $P_3$.   Noting that $-\Delta=D^2$, we have 
$$ P_3=\langle (-\Delta \rho) \  | \ (\nabla \cdot \vec{X})\rho\rangle=\langle D^2\rho \  | \ D^{-1}\rho\rangle=
\langle D^{\frac{1}{2}}\rho \  | \ D^{\frac{1}{2}}\rho\rangle=\frac{1}{2}\|D^{\frac{1}{2}}(|u|^2)\|_{L^2}^2 \, .$$ 
There are two terms for which their positivity is not immediate--the terms $P_{2}$ and $P_{4}$.  To check $P_2$, recall \eqref{E:dX}, where the kernel $\eta_{kj}(x,y)$ is symmetric.  Then
$$P_2=\int_{\Bbb R^2 \times \Bbb R^2}\left( \frac{\rho(y)}{\rho(x)}p_k(x)p^j(x)-p_k(y)p^j(x) \right) \eta_{j}^{k}(x,y)\, dxdy \, .$$
By symmetry,  we get
$$P_2=\int_{\Bbb R^2 \times \Bbb R^2}\left( \frac{\rho(x)}{\rho(y)}p_k(y)p^j(y)-p_k(x)p^j(y) \right)\eta_{j}^{k}(x,y) \, dxdy \, .$$
and thus,
\begin{align*}
P_2 &= \frac{1}{2}\int_{\Bbb R^2 \times \Bbb R^2}\left( \frac{\rho(y)}{\rho(x)}p_k(x)p^j(x)+\frac{\rho(x)}{\rho(y)}p_k(y)p^j(y)-p^j(x)p_k(y)-p^j(y)p_k(x) \right) \eta_{j}^{k}(x,y) \, dxdy \\
& =\frac{1}{2}\int_{\Bbb R^n \times \Bbb R^n}\left( \sqrt{\frac{\rho(y)}{\rho(x)}}p_k(x)-\sqrt{\frac{\rho(x)}{\rho(y)}}p_k(y)\right)
\left( \sqrt{\frac{\rho(y)}{\rho(x)}}p^j(x)-\sqrt{\frac{\rho(x)}{\rho(y)}}p^j(y)\right) \eta_{j}^{k}(x,y) \, dxdy.
\end{align*}
Thus,  if we define the two point momentum vector
$$\vec{J}(x,y)=\sqrt{\frac{\rho(y)}{\rho(x)}}\vec{p}(x)-\sqrt{\frac{\rho(x)}{\rho(y)}}\vec{p}(y) \, ,$$
we can write
$$P_2=\frac{1}{2}\langle J^{j}J_{k} \  | \ (\partial_{j}X^{k})\rangle \geq 0 \, ,$$
since $\partial_{j}X^{k}$ is positive definite. 

Now we discuss $P_4$.  Recalling the definition of $V$ and $\vec{X}$, we expand the form and obtain
$$P_4=4\gamma \int_{\Bbb R^2}\int_{\Bbb R^2}\int_{\Bbb R^2}\rho(x)\rho(y)\rho(z)\frac{(x-y)\cdot (x-z)}{|x-y|\, |x-z|^{\gamma +2}} \, dxdydz$$
By symmetry, this term becomes
$$P_4=2\gamma \int_{\Bbb R^2}\int_{\Bbb R^2}\int_{\Bbb R^2}\frac{\rho(x)\rho(y)\rho(z)}{|x-z|^{\gamma +2}}\left( \frac{(x-y)\cdot (x-z)}{|x-y|}+\frac{(z-y)\cdot (z-x)}{|z-y|}\right) \, dxdydz \, .$$ 
But
\begin{align*}
\indentalign  \frac{(x-y)\cdot (x-z)}{|x-y|}+\frac{(z-y)\cdot (z-x)}{|z-y|} \\
&=|x-y|+\frac{(x-y)\cdot (y-z)}{|x-y|}+|z-y|+\frac{(y-x)\cdot (z-y)}{|z-y|} 
\end{align*}
Since the 2nd term is dominated by the 3rd and the 4th is dominated by the 1st, the sum is greater than $0$.  Thus, $P_4\geq 0$. 

Since all the terms are positive, we keep only $P_3$, and after integrating in time we have 
$$\|D^{\frac{1}{2}}(|u|^2)\|_{L_{t}^{2}L_{x}^2}^2 \lesssim \sup_{t} M(t).$$
It remains to show that $M(t)$ is bounded by the appropriate norms.  But
\begin{align*}
\frac{1}{2}M(t) &= \langle [x;D^{-1}]^{j}\rho(t) \  | \ p_{j}(t)\rangle \lesssim \|p_{j}\|_{L^1}\|[x;D^{-1}]_j\rho(t)\|_{L^{\infty}} \\
&\lesssim \|p_{j}\|_{L^1}\|\rho\|_{L^1}\|[x;D^{-1}]_j\|_{L^1\rightarrow L^{\infty}} \, .
\end{align*}
We clearly have that  $\|p_{j}\|_{L^1}\lesssim \|u\|_{\dot{H}^{\frac{1}{2}}}^2$ and $\|\rho\|_{L^1}=\frac{1}{2}\|u\|_{L^{2}}^2$.  Finally, the operator norm $\|[x;D^{-1}]_j\|_{L^1\rightarrow L^{\infty}}$ is bounded by $1$, which is immediate from the definition \eqref{E:Xint} of $X$.

Thus, all in all, we have \eqref{22d}.
\end{proof}

In this paper, for both the nonlinear Schr\"odinger and Hartree equations, we will use the following estimate, which is an easy consequence of Theorem \ref{theorem3} and interpolation.

\begin{corollary}[Morawetz interpolated space-time bounds for NLS and Hartree] 
\label{C:Morawetz_interpolated}
Let $u$ be a solution to NLS with $n=2$, $p>1$ or Hartree with $n=2$, $0<\gamma<2$.  Then for
$$ \frac3q+\frac2r=1, \quad 4\leq q \leq \infty, \quad 2\leq r \leq 8\,, $$
we have the \emph{a priori} bound
$$\|u\|_{L_t^qL_x^r} \lesssim E(u_{0})^{\frac{r-2}{6r}}\|u_{0}\|_{L^2}^{\frac{2(r+1)}{3r}}$$
\end{corollary}

\begin{proof}
By Sobolev imbedding and Theorem \ref{theorem3}, we have
$$\| u\|_{L_t^4L_x^8}^2 = \| |u|^2 \|_{L_t^2L_x^4} \lesssim \|D^{1/2} |u|^2 \|_{L_t^2L_x^2} \lesssim \|u\|_{L_t^\infty \dot H_x^{1/2}} \|u\|_{L_t^\infty \dot L_x^2} \lesssim E(u_{0})^{\frac{1}{4}}\|u_{0}\|_{L^2}^{\frac{3}{2}}$$
We conclude by using the interpolation (H\"older) inequality
$$\|u\|_{L_t^qL_x^r} \leq \|u\|_{L_t^\infty L_x^2}^{1-\theta} \|u\|_{L_t^4L_x^8}^\theta \,, \quad \theta = \frac43-\frac8{3r}.$$
\end{proof}

\section{Weak time-decay bounds}

In this section we establish the key estimates required to implement Nakanishi's argument.  These estimates are consequences of the Morawetz bounds derived in \S\ref{S:Morawetz} (in particular, Corollary \ref{C:Morawetz_interpolated}) and the $L_x^p$ time-decay estimates for the linear group $U(t)$.

\begin{lemma}[Weak time-decay bound for NLS]
\label{L:NLSwkdecay}
Suppose $n=2$, $2<p<\infty$, and $u(t)$ is a (global) $H^1$ solution to the defocusing NLS.  Then there exists $\alpha,\beta>0$ such that for any $\psi \in C_c^\infty$ and any $s<t$, we have
$$\left| \int_s^t \la  |u(\sigma)|^{p-1}u(\sigma), U(\sigma)\psi \ra \, d\sigma \right| \lesssim \begin{cases}  |t-s|^\alpha & \\ s^{-\beta} & \text{if }s\geq 1 \end{cases}$$
This notation means that both estimates hold.  The brackets $\la \cdot, \cdot \ra$ represent the $L^2_x$ spatial pairing.  The implicit constant in this estimate depends on $\|\psi\|_{L^1}$, $\|\psi\|_{L^2}$, and $\|u\|_{L_{t\in  (-\infty,+\infty)}^\infty H_x^1}$.  We recall that $\|u\|_{L_{t\in  (-\infty,+\infty)}^\infty H_x^1}$ is controlled by the mass and energy.
\end{lemma}
\begin{proof}
We begin by proving the first estimate.  By the unitarity $\|U(\sigma)\psi\|_{L^2} = \|\psi\|_{L^2}$, we have
$$\left| \int_s^t \la  |u(\sigma)|^{p-1}u(\sigma), U(\sigma)\psi \ra \, d\sigma \right| \lesssim \int_s^t \|u(\sigma)\|_{L^{2p}}^p \, d\sigma$$
Employing the Gagliardo-Nirenberg estimate $\|u(\sigma)\|_{L_x^{2p}}^p \leq \|u\|_{L^2} \|\nabla u\|_{L^2}^{p-1}$, we obtain that the above display is controlled by $|s-t|$.

Now we derive the second estimate.  Suppose $s\geq 1$.  By the 2-d space-time decay estimate $\|U(\sigma)\psi\|_{L_x^\infty} \leq \sigma^{-1}\|\psi\|_{L^1}$, we obtain
$$\left| \int_s^t \la  |u(\sigma)|^{p-1}u(\sigma), U(\sigma) \psi \ra \, d\sigma \right| \lesssim  \int_s^t \sigma^{-1} \|u(\sigma)\|_{L_x^p}^p\, d\sigma$$
Let $q$ be such that $\frac3q=1-\frac2p$.  First, suppose $p<5$.  Then $q>p$, and we can apply H\"older in $\sigma$ with conjugate pair $(\frac{3}{5-p},\frac{q}{p})$ to obtain that the above display is controlled by 
$s^{-\frac{p-2}3}\|u\|_{L_{[s,t]}^qL_x^p}^p$.  Corollary \ref{C:Morawetz_interpolated} furnishes the required bound.  Suppose now that $p\geq 5$.  Then, by the Gagliardo-Nirenberg inequality, 
$\|u(\sigma)\|_{L_x^p}^p \leq \|u(\sigma)\|_{L_x^5}^5\|\nabla u(\sigma)\|_{L_x^2}^{p-5}$, and therefore
$$\int_s^t \sigma^{-1}  \|u(\sigma)\|_{L_x^p}^p\, d\sigma \leq s^{-1}\|\nabla u\|_{L_t^\infty L_x^2}^{p-5}\|u\|_{L_t^5L_x^5}^5 \,.$$
We now appeal to Corollary \ref{C:Morawetz_interpolated} with $r=q=5$.

\end{proof}

\begin{lemma}[Weak time-decay bound for Hartree]
\label{L:Hwkdecay}
Suppose $n=2$, $1<\gamma <2$, and $u(t)$ is a (global) $H^1$ solution to the Hartree equation.  Then there exists $\alpha,\beta>0$ such that for any $\psi \in C_c^\infty$ and any $s<t$, we have
$$\left| \int_s^t \la  (|x|^{-\gamma}\star|u(\sigma)|^2)u(\sigma), U(\sigma)\psi \ra \, d\sigma \right| \lesssim \begin{cases}  |t-s|^\alpha & \\ s^{-\beta} & \text{if }s\geq 1 \end{cases}$$
This notation means that both estimates hold.  The brackets $\la \cdot, \cdot \ra$ represent the $L^2_x$ spatial pairing.  The implicit constant in this estimate depends on $\|\psi\|_{L^1}$, $\|\psi\|_{L^\frac43}$, and $\|u\|_{L_{t\in  (-\infty,+\infty)}^\infty H_x^1}$.  We recall that $\|u\|_{L_{t\in  (-\infty,+\infty)}^\infty H_x^1}$ is controlled by the mass and energy.
\end{lemma}
\begin{proof}
By the Hardy-Littlewood-Sobolev (HLS) inequality, \cite{es}, we have that
\begin{equation}
\label{E:HLS}
\| |x|^{-\gamma}\star |u|^2 \|_{L^2} \lesssim \| |u|^2 \|_{L^\frac{2}{3-\gamma}} = \|u\|_{L^\frac{4}{3-\gamma}}^2
\end{equation}
HLS in 2-d requires $\gamma<2$, and the restriction to $\gamma>1$ implies that $\frac{2}{3-\gamma}>1$, which is also required for the validity of HLS.

To derive the first estimate, we apply H\"older with partition $(2,4,4)$, and use the space-time estimate $\|U(\sigma)\psi\|_{L^4} \lesssim |\sigma|^{-1/2} \|\psi\|_{L^{4/3}}$ to obtain
$$\left| \int_s^t \la  (|x|^{-\gamma}\star|u(\sigma)|^2)u(\sigma), U(\sigma)\psi \ra \, d\sigma \right| \lesssim \int_s^t \|(|x|^{-\gamma}\star|u(\sigma)|^2)\|_{L_x^2}\|u(\sigma)\|_{L_x^4} |\sigma|^{-1/2} \,d \sigma $$
We now follow through with \eqref{E:HLS} and Gagliardo-Nirenberg to obtain the bound $\int_s^t |\sigma|^{-1/2} \, d\sigma$.  If $0\leq s<t$, this evaluates to $t^{1/2}-s^{1/2} \leq (t-s)^{1/2}$.  If $s\leq t\leq 0$, then we similarly obtain a bound of $|t-s|^\frac12$.  If $s<0<t$, then the integral evaluates to $|s|^\frac12+t^\frac12 \leq 2|t-s|^\frac12$.

To derive the second estimate, we apply H\"older with partition $(2,2,\infty)$,  and use the space-time estimate $\|U(\sigma)\psi\|_{L^\infty} \lesssim |\sigma|^{-1} \|\psi\|_{L^1}$ to obtain 
$$\left| \int_s^t \la  (|x|^{-\gamma}\star|u(\sigma)|^2)u(\sigma), U(\sigma)\psi \ra \, d\sigma \right| \lesssim \int_s^t \|(|x|^{-\gamma}\star|u(\sigma)|^2)\|_{L_x^2}\|u(\sigma)\|_{L_x^2} |\sigma|^{-1} \,d \sigma $$
Following through with \eqref{E:HLS}, we obtain
$$\lesssim \int_s^t |\sigma|^{-1}\|u(\sigma)\|_{L_x^\frac{4}{3-\gamma}}^2 \, d\sigma$$
Set $r=\frac{4}{3-\gamma}$, note that ($1<\gamma<2 \implies 2<r<4$) and define $q$ so that $\frac{3}{q}+\frac2{r}=1$.  Since $q>2$ (in fact $q>6$) we can apply H\"older in $\sigma$ with partition $(\frac{3}{4-\gamma}, \frac{q}{2})$ to obtain that the above is bounded by $s^{-\frac{\gamma-1}{3}}\|u\|_{L_t^qL_x^r}^2$.  Corollary \ref{C:Morawetz_interpolated} furnishes the required bound.
\end{proof}

\section{Basic uniqueness and regularity}

A key step in the proof of the main theorem is the statement that if $u$ solves NLS or Hartree as a distribution with some additional regularity, then $u$ is in fact the unique (strong) solution satisfying the mass and energy conservation laws.  In this section we record precise statements.   While all the material in this section appears in the literature (see e.g. Prop. 4.2.9, Lemma 4.2.8, and Theorem 3.3.9 in \cite{tc} and Theorems I, III in Kato \cite{tk}), we believe that the $H^1$ \emph{unconditional} uniqueness statements are less well-known than the $H^1$ local well-posedness statements that come equipped with \emph{conditional} uniqueness.  Therefore, we will effectively assume the standard $H^1$ local existence and carefully discuss the unconditional uniqueness claims.

We will need the Strichartz estimates in the following form.  Suppose that $(r,q)$ satisfy the admissibility condition
$$\frac{2}{q}+\frac{n}{r}=\frac{n}{2}, \qquad 2\leq q\leq \infty$$
excluding the case $n=2$ and $q=2$.  Then
$$\|e^{it\Delta}\phi \|_{L_t^qL_x^r} \lesssim \|\phi\|_{L^2}$$
$$\left\|\int_0^t e^{i(t-t')\Delta}f(\cdot,t') \,dt' \right\|_{L_t^qL_x^r} \lesssim \|f\|_{L_t^{q'}L_x^{r'}}$$
where $q'$, $r'$ are the H\"older duals to $q$, $r$, respectively.

Let $I\subset \mathbb{R}$ be a bounded\footnote{For convenience in the discussion that follows, we require $I$ to be bounded.  Of course, one can concatenate intervals to draw global in time conclusions.} open interval of time.  Since we are discussing distributional solutions to NLS and Hartree with $u\in L_I^\infty H_x^1$, we first off need to make sure that the nonlinearities belong to $L_{\text{loc}}^1$ (and are thus well-defined distributions) and that the interaction potential term in the energy is well-defined.  

For NLS, we assume that $p$ is energy-subcritical, i.e. 
\begin{equation}
\label{E:energysub}
1<p<\frac{n+2}{n-2}
\end{equation}
Then by Sobolev imbedding, $u\in L_I^\infty L_x^r$ for all $2\leq r < \frac{2n}{n-2}$.  Thus, the nonlinearity $|u|^{p-1}u\in L_{\text{loc}}^1$, and furthermore $u\in L_I^\infty L_x^{p+1}$, which implies that the energy is a well-defined quantity for a.e. $t$.  

For Hartree, suppose $0<\gamma<n$.  Note that by the Hardy-Littlewood-Sobolev inequality, if $u(t)\in H_x^1$, then $|x|^{-\gamma}\star |u(t)|^2 \in L_x^p$ for all $p$ such that
$$\max\left( 1-\frac2n-\frac{n-\gamma}{n}, 0\right) < \frac1p < 1-\frac{n-\gamma}{n}$$
Thus, we certainly have that the nonlinearity $(|x|^{-\gamma}\star|u|^2)u \in L_{\text{loc}}^1$ is a well-defined distribution.

\begin{lemma}
\label{L:interp}
Suppose that $I$ is an open interval of time, $u\in L_I^\infty H_x^1$ and $\partial_t u \in L_I^\infty H_x^{-1}$. (We are not asserting here that $u$ solves any equation.)  Then $u$ can be redefined on a set of zero measure in time such that $u \in C_I^{0,\frac12-\frac{s}{2}}(I; H_x^s)$ for all $-1\leq s<1$, with
$$\|u(t_2)-u(t_1)\|_{H_x^s} \lesssim |t_2-t_1|^{\frac12-\frac{s}{2}}\|\partial_t u\|_{L_I^\infty H_x^{-1}}^{\frac12-\frac{s}{2}}\|u\|_{L_I^\infty H_x^1}^{\frac12+\frac{s}{2}} \,.$$
Note that we do not claim that it follows that $u\in C(I;H_x^1)$ under these hypotheses alone (the case $s=1$ is excluded), nor do we claim that $u\in C^1(I;H_x^{-1})$.  It is true, however, that $u\in C^{0,1}(I;H_x^{-1})$, i.e. $u$ is Lipschitz continuous as a map from time into $H_x^{-1}$ (the case $s=-1$ is included in the above estimate).
\end{lemma}

\begin{proof}
For $s=-1$ the result is immediate by the fundamental theorem of calculus and for $s=1$, the statement is trivial.  The result follows by interpolation.
\end{proof}

\begin{lemma}[Uniqueness and regularity for NLS]
\label{L:NLS_uniqueness} 
Suppose $u$ solves NLS (in the sense of distributions) with $1<p<\frac{n+2}{n-2}$ on a bounded open time interval $I$ and $u\in L_I^\infty H_x^1$. Then $u(t)$ can be redefined on a set of times $t$ of zero measure such that $u\in C(I;H_x^1)\cap C^1(I;H_x^{-1})\cap L_I^qW_x^{1,p+1}$, and $u$ satisfies the mass and energy conservation laws.  Here, $2<q\leq \infty$ is defined by
$$\frac{2}{q}+\frac{n}{p+1}=\frac{n}{2} \,.$$
Also, if $v$ is another distributional solution to NLS on $I$, $v\in L_I^\infty H_x^1$,  and $v(t)$ and $u(t)$ agree at some $t_0\in I$ (after being redefined as above so that $v,u\in C(I;H_x^1)$), then $u\equiv v$.
\end{lemma}

\begin{proof}
Clearly, $\Delta u \in L_I^\infty H_x^{-1}$.  By Sobolev imbedding, 
$$\| |u|^{p-1}u \|_{H_x^{-1}} \lesssim \| |u|^{p-1}u \|_{L_x^\alpha} = \|u\|_{L_x^{\alpha p}}^p \,$$
for any $\alpha$ such that $\max( \frac{2n}{n+2}, \,1+\, ) \leq \alpha \leq 2$.   Now if $2\leq \alpha p \leq \frac{2n}{n-2}$, then  by Gagliardo-Nirenberg-Sobolev, we'll have $\|u\|_{L_x^{\alpha p}} \lesssim \|u \|_{H_x^1}$.  But it is clear that for $1<p<\frac{n+2}{n-2}$, the two intervals $\max( \frac{2n}{n+2}, \, 1+\, ) \leq \alpha \leq 2$ and $\frac{2}{p} \leq \alpha \leq \frac{2n}{(n-2)p}$ have a nontrivial intersection.

Consequently, $\partial_t u = -\Delta u +|u|^{p-1}u \in L_I^\infty H_x^{-1}$, and we can apply Lemma \ref{L:interp} to $u$ to conclude that $u\in C^{0,\frac12-\frac{s}{2}}(I; H_x^s)$ for all $-1\leq s <1$.  In particular, $u(t)$ is now well-defined as an element of $L_x^2$ for \emph{all} $t$.  Moreover, since $u\in L_I^\infty H_x^1$, we have that $u(t)\in H_x^1$ for all $t$ except possibly on a set of measure zero.  Select some $t_0\in I$ such that $u(t)\in H_x^1$.  By the standard local theory (see Kato \cite{tk}, Theorems I, III), there exists an interval $I'$ containing $t_0$ with length depending only on $\|u(t_0)\|_{H_x^1}$, and a solution $v\in C^1(I';H_x^{-1})\cap C(I';H_x^1) \cap L_{I'}^qW_x^{1,p+1}$ of NLS such that $v(t_0)=u(t_0)$ and $\|v\|_{C(I';H_x^1)} \lesssim \|u(t_0)\|_{H_x^1}$.  Kato \cite{tk} also shows the solution $v$ satisfies the mass and energy conservation laws.  By appealing to the claim below (with $I$ replaced by $I'$), we conclude that $u=v$ on $I'$.  This completes the proof, since we can carry out this argument on subintervals $I'$ that fill out $I$.

\smallskip

\noindent\emph{Claim}. If $u$ and $v$ have the properties given in the theorem statement and there exists $t_0\in I$ such that $v(t_0)=u(t_0)$, then necessarily $u=v$ on $I$. (Here, the statement $u(t_0)=v(t_0)$ is well-defined since we know by the above argument that both $u,v\in C(I; L_x^2)$).  

\smallskip

To prove the claim, let $w=u-v$.  We will show that on some interval $I'$ containing $t_0$ with $|I'|$ depending only on $\max(\|u\|_{L_I^\infty H_x^1}, \|v\|_{L_I^\infty H_x^1})$, we have $w\equiv 0$.  The argument can then be iterated  to obtain that $u=v$ on all of $I$.  Since $U(t)$ is a strongly continuous unitary group on $H_x^{-1}$, $w(t)$ satisfies the integral equation
$$w(t) =  i \int_{t_0}^t U(t'-t_0) (|u(t')|^{p-1}u(t')-|v(t')|^{p-1}v(t')) \, dt' \,.$$
Applying the Strichartz estimates with $r=p+1$,
\begin{align*}
 \|w\|_{L_{I'}^qL_x^{p+1}} 
&\lesssim \| |u|^{p-1} u - |v|^{p-1}v \|_{L_{I'}^{q'}L_x^{\frac{p+1}{p}}} \\
& \lesssim |I'|^{\frac{1}{q'}-\frac1{q}}(\|u\|_{L_{I'}^\infty L_x^{p+1}}^{p-1}+\|v\|_{L_{I'}^\infty L_x^{p+1}}^{p-1}) \|w\|_{L_{I'}^qL_x^{p+1}}
\end{align*}
Taking $|I'|$ small depending only on $\max(\|u\|_{L_I^\infty H_x^1}, \|v\|_{L_I^\infty H_X^1})$, we obtain that $w\equiv 0$ on $I'$, as claimed.

\end{proof}

\begin{lemma}[Uniqueness and regularity for Hartree]
\label{L:Hartree_uniqueness} 
Suppose $u$ solves Hartree (in the sense of distributions) with $0<\gamma<\min(n,4)$ on a bounded open time interval $I$ and $u\in L_I^\infty H_x^1$.  Then $u(t)$ can be redefined on a set of times $t$ of zero measure such that $u\in C(I;H_x^1)\cap C^1(I;H_x^{-1})\cap L_I^qW_x^{1,r}$, and $u$ satisfies the mass and energy conservation laws.  Here, $r$ and $q$ are defined by
$$\frac1r = \frac12-\frac{\gamma}{4n}, \qquad \frac{2}{q}+\frac{n}{r}=\frac{n}{2} \,.$$
Also, if $v$ is another distributional solution to Hartree on $I$, $v\in L_I^\infty H_x^1$,  and $v(t)$ and $u(t)$ agree at some $t_0\in I$ (after being redefined as above so that $u,v\in C(I;H_x^1)$), then $u\equiv v$.
\end{lemma}
\begin{proof}
By Sobolev imbedding, $\| (|x|^{-\gamma}\ast |u|^2)u\|_{H_x^{-1}} \leq \| (|x|^{-\gamma}\ast |u|^2)u\|_{L_x^\alpha}$ for any $\alpha$ such that 
$$\frac12 \leq \frac1{\alpha} \leq \min\Big( \frac12+\frac1n \, , \; 1-\Big)$$
Select an $\alpha$ meeting this condition and, in addition, the following:
$$\max\Big(\frac12+\frac{\gamma}{n}-\frac{3}{n} \, , \; \frac12-\frac{\gamma}{2n}+ \, \Big) \leq \frac1{\alpha}\leq \frac{1}{2}+\frac{\gamma}{n} .$$
This is possible since $0<\gamma< \min(\gamma, 4)$.  Let
$$\frac1{\beta} = \frac13\Big(1-\frac{\gamma}{n}+\frac1{\alpha}\Big), \qquad \frac{1}{\tilde \beta} = \frac{2}{\beta}-1+\frac{\gamma}{n} \,.$$ 
We note that $2\leq \beta \leq \frac{2n}{n-2}$, $\tilde \beta <\infty$, and $\frac{1}{\tilde \beta} + \frac{1}{\beta} = \frac{1}{\alpha}$.  Applying H\"older and the Hardy-Littlewood-Sobolev theorem on fractional integration, we have
$$\| (|x|^{-\gamma}\ast |u|^2)u\|_{L_x^\alpha} \leq \|(|x|^{-\gamma}\ast |u|^2)\|_{L_x^{\tilde \beta}} \|u\|_{L_x^\beta} \lesssim \|u\|_{L_x^\beta}^3 \,.$$
By Sobolev imbedding, we have $\|u\|_{L_x^\beta} \lesssim \|u\|_{H_x^1}$.  Thus, $(|x|^{-\gamma}\ast |u|^2) u \in L_I^\infty H_x^{-1}$.  It is also clear that $\Delta u \in L_I^\infty H_x^{-1}$, and therefore from the Hartree equation it follows that $\partial_t u \in L_I^\infty H_x^{-1}$.  We can thus appeal to Lemma \ref{L:interp}, and conclude that $u$ can be redefined on a set of times $t$ of zero measure such that $u\in C^{0,\frac12-\frac{s}{2}}(I;H_x^1)$ for $-1\leq s<1$.  In particular, $u(t)$ is well-defined for \emph{all} $t$ as an element of $L_x^2$.  Moreover, since $u\in L_I^\infty H_x^{-1}$, we have, except on a set of measure zero, that $u(t)\in H_x^1$.  Take any $t_0\in I$ for which $u(t)\in H_x^1$.  By the local theory for the Hartree equation (see Ginibre-Velo \cite{gv-hartree-lwp}), there exists an interval $I'$, with center $t_0$ and $|I'|$ depending only on $\|u(t_0)\|_{H_x^1}$, and  $v\in C^1(I'; H_x^{-1})\cap C(I'; H_x^1) \cap L_{I'}^{q}W_x^{1,r}$ solving the Hartree equation, with $v(t_0)=u(t_0)$.  Moreover, $v$ satisfies the energy and mass conservation laws.  The proof is completed upon appealing to the following uniqueness claim.

\smallskip

\noindent \emph{Claim}:  Suppose that $u$, $v$ are given functions satisfying the conditions in the theorem statement and there exists $t_0\in I$ such that $u(t_0)=v(t_0)$.  (The statement $u(t_0)=v(t_0)$ is well-defined since the above argument established that $u,v\in C(I;L_x^2)$).  Then $u\equiv v$ on $I$.

\smallskip

Let $w=u-v$.  By iteration, it suffices to show that $w\equiv 0$ on an interval $I'$ centered at $t_0$ with length depending only on $\max(\|u\|_{L_I^\infty H_x^1}, \|v\|_{L_I^\infty H_x^1})$.  

We note that with $r$ as defined in the statement of the theorem, 
$$\max\Big(\frac14,\frac12-\frac1n\Big) < \frac1r<\frac12 \,,$$
and so $(r,q)$ is a Strichartz admissible pair.  Let $\tilde r$ be defined by
$$\frac1{\tilde r} = \frac2r - \frac{n-\gamma}{n}=\frac{\gamma}{2n}\, .$$
By the Hardy-Littlewood-Sobolev inequality,
\begin{equation}
\label{E:uniqueHartree1}
\| |x|^{-\gamma} \ast (u_1 u_2) \|_{L_x^{\tilde r}} \lesssim \| u_1u_2 \|_{L_x^{\frac{r}{2}}} \leq \| u_1\|_{L_x^r}\|u_2\|_{L_x^r}
\end{equation}
Moreover, $\frac1{r}+\frac1{\tilde r} = \frac1{r'}$. (Recall $r'$ is the H\"older dual to $r$.)  Consequently, by H\"older and \eqref{E:uniqueHartree1} we have
\begin{equation}
\label{E:uniqueHartree2}
\| (|x|^{-\gamma}*(u_1\bar u_2))u_3 \|_{L_x^{r'}} \lesssim  \|u_1\|_{L_x^r} \|u_2\|_{L_x^r} \|u_3\|_{L_x^r}
\end{equation}
Since $U(t)$ is a strongly continuous unitary group on $H_x^{-1}$, we have that $w$ satisfies the associated integral equation.  Applying the Strichartz estimates, we obtain
$$\|w\|_{L_{I'}^q L_x^r} \lesssim \| (|x|^{-\gamma}*|u|^2)u - (|x|^{-\gamma}*|v|^2)v \|_{L_{I'}^{q'}L_x^{r'}}$$
Applying \eqref{E:uniqueHartree2}, 
$$\|w\|_{L_{I'}^q L_x^r} \lesssim |I|^{\frac{1}{q'}-\frac{1}{q}} (\| u\|_{L_{I'}^\infty L_x^r}^2 + \|v\|_{L_{I'}^\infty L_x^r}^2)\|w\|_{L_{I'}^qL_x^r} \,.$$
Select $|I'|$ sufficiently small, depending only on $\max( \|u\|_{L_I^\infty H_x^1}, \|v\|_{L_I^\infty H_x^1})$, to obtain $w\equiv 0$ on $I'$.
\end{proof}

\section{Nakanishi's argument}

We will prove Theorem \ref{Theorem1}.   The only difference between the NLS case and the Hartree case is that, for the former, we use Lemma \ref{L:NLSwkdecay}, \ref{L:NLS_uniqueness} and for the latter, we use Lemma \ref{L:Hwkdecay}, \ref{L:Hartree_uniqueness}.  Therefore, we shall only present the details for the NLS case.  

We are given an $H_x^1$ solution $v(t)$ to the linear Schr\"odinger equation and want to find an $H_x^1$ solution $u(t)$ to NLS such that $\|u(t)-v(t)\|_{H_x^1} \to 0$ as $t\to +\infty$.  For any time $T>0$, we have a unique solution $w_T=w \in C(\Bbb R ; H^{1})$ to the initial value problem
\begin{equation}
\left\{
\begin{matrix}
iw_{t}+ \Delta w -|w|^{p-1}w=0, & x \in {\mathbb R^2}, & t\in {\mathbb R},\\
w(x,T)=v(T)\in H^{1}({\mathbb R^2}).
\end{matrix}
\right.
\end{equation} 
We start by constructing $u(t)$ as a weak limit of $w_T(t)$ as $T\to +\infty$ and establish that $U(-t)u(t)$ converges weakly to $v(0)$ as $t\to \infty$ in $H_x^1$.  This is accomplished by an argument using the Arzela-Ascoli theorem and the (uniform in $n$) weak decay estimate in Lemma \ref{L:NLSwkdecay}.  Then, by using the conservation laws, we upgrade the weak convergences to strong convergences, and obtain in particular that $U(-t)u(t) \to v(0)$ strongly in $H_x^1$ as $t\to \infty$.

We now turn to the details. Let $\{ \psi_j \}_{j=1}^\infty \subset C_c^\infty(\mathbb{R}^n)$ be a basis of $H^{-1}(\mathbb{R}^n)$ and let $f_{j,T}(t) = \la U(-t)w_T(t),\psi_j \ra$.  

\medskip

\noindent\textbf{Assertion 1}.   $\forall \, j\in \mathbb{N}$, $\lim_{t\to \pm \infty} f_{j,T}(t)$ exists (and is a number in $\mathbb{C}$).  Thus, $f_{j,T}$ is defined on the compact interval $[-\infty,+\infty]$.

\smallskip

\noindent \textit{Proof of Assertion 1}.  By writing the equivalent integral representation of the solution $w_{T}(t)$ we have that
$$w_{T}(t)=U(t-T)v(T)-i\int_{T}^{t}U(t-\sigma)\{|w_{T}|^{p-1}w_{T}(\sigma)\}d\sigma$$ 
where $U(t)=e^{it\Delta}$ is the linear semigroup. Then
\begin{align*}
U(-t)w_{T}(t) & =U(-T)v(T)-i\int_{T}^{t}U(-\sigma)\{|w_{T}|^{p-1}w_{T}(\sigma)\}d\sigma \\
&= v(0)-i\int_{T}^{t}U(-\sigma)\{|w_{T}|^{p-1}w_{T}(\sigma)\}d\sigma
\end{align*}
and thus
$$\langle U(-t)w_{T}(t), \psi_{j}\rangle=\langle v(0), \psi_{j}\rangle-i\int_{T}^{t}\langle U(-\sigma)\{|w_{T}|^{p-1}w_{T}(\sigma)\}, \psi_{j} \rangle d\sigma.$$
Now for any $T \geq 1$, by the proof of Lemma \ref{L:NLSwkdecay}, we have that the limit of $f_{j,T}(t) = \la U(-t)w_T(t),\psi_j \ra$ exists as $t \rightarrow \pm \infty$ and it is a finite number.  \textit{End proof of Assertion 1}.
 
\medskip

\noindent\textbf{Assertion 2}.  $\forall \, j\in \mathbb{N}$, the family $\{f_{j,T}\}_T$ is equicontinuous in $C([-\infty,+\infty])$ (here $C([-\infty,+\infty])$ is given the sup norm over $t\in[-\infty,+\infty]$).  Moreover, $\forall \, j\in \mathbb{N}$, the family $\{f_{j,T}\}_T$ is uniformly bounded.

\smallskip

\noindent\textit{Proof of Assertion 2}. The equicontinuity of $\{f_{j,T}\}_T$ follows immediately from Lemma 3.1. Moreover we have that
\begin{align*}
|f_{j,T}(t)| = |\la & U(-t)w_T(t),\psi_j \ra| \lesssim \|w_{T}(t)\|_{H^1} \lesssim E^{\frac{1}{2}}(w_{T})(t)  \\
& =E^{\frac{1}{2}}(w_T)(T) = E^{\frac12}(v)(T) \\
&  = \left ( \frac{1}{2}\int |\nabla v(T)|^{2}dx+\frac{1}{p+1}\int |v(T)|^{p+1}dx \right )^{\frac{1}{2}} \\
&\lesssim \|v(0)\|_{H^1} + \|v(0)\|_{H^1}^\frac{p+1}{2} \,.
\end{align*}
The last inequality follows because $\|v(T)\|_{L^2}=\|v(0)\|_{L^2}$, $\|\nabla v(T)\|_{L^2} = \|\nabla v(0)\|_{L^2}$, and by Gagliardo-Nirenberg, we have that
$$\|v(T)\|_{L^{p+1}}\lesssim \|v(T)\|_{L^2}^{\frac{2}{p+1}}\|\nabla v(t)\|_{L^{2}}^{\frac{p-1}{p+1}} \lesssim \|v(0)\|_{H^{1}}.$$
Thus $\{f_{j,T}\}_T$ is uniformly in $T$, bounded. \textit{End proof of Assertion 2}.

\medskip

Now, by Assertions 1--2 and the Arzela-Ascoli theorem and a diagonal argument, there exists a sequence $T_n\to +\infty$ such that $\forall \, j\in \mathbb{N}$, the sequence $f_{j,T_n}$ converges strongly in $C([-\infty,+\infty])$ as $n\to \infty$ (uniformly in $t$).  This is the key ingredient in the proof of Assertions 3 and 6 below.

\medskip

\noindent\textbf{Assertion 3}.  For each $t\in (-\infty,+\infty)$, the sequence $w_{T_n}(t)$ converges weakly in $H^1$ as $n \rightarrow \infty$.  Let $u(t)$ denote this weak limit. 

\medskip

It is important to note one strength of Assertion 3:  there is a \emph{single} sequence $T_n$ such that \emph{for all} $t$, the weak convergence along this sequence $T_n$ holds.  In contrast, it is a simple statement that for each $t$, there is a sequence $T_n$ for which the weak convergence holds, with the sequence $T_n$ depending on $t$.

\medskip

\noindent\textit{Proof of Assertion 3}.  \underline{Fix any $t$}.   We know, by the above Arzela-Ascoli argument, that $\forall \; j$, $\la w_{T_{n}}(t), U(t)\psi_j \ra$ converges as $n \rightarrow \infty$ to a limit; let's call it $L(U(t)\psi_{j})$.  $L$ clearly extends to a linear functional on the finite span $E\subset H^{-1}$ of $\{ U(t)\psi_j \}_j$, and $L: E\to \mathbb{C}$ is bounded since $w_{T_{n}}(t)$ is bounded in $H^{1}$ uniformly in $T_{n}$.  Since $E$ is dense in $H^{-1}$, $L$ extends to a bounded linear functional on all of $H^{-1}$.  Then by the Riesz Representation Theorem we know that there exists some $u(t) \in H^1$ such that $L(U(t)\psi_{j})=\la  u(t), U(t)\psi_{j}\ra$. Thus
$$\lim_{n \rightarrow \infty}\la U(-t)w_{T_{n}}(t),\psi_j \ra=\la  U(-t)u(t), \psi_{j}\ra \,,$$
i.e. $U(-t)w_{T_n}(t)$ converges weakly in $H^1$ to $U(-t)u(t)$ as $n\to +\infty$.  Since $U(t)$ is unitary on $H^{-1}$, $w_{T_n}(t)$ converges weakly in $H^1$ to  $u(t)$ as $n\to +\infty$.  \textit{End proof of Assertion 3}.

\medskip

The quantities $\|w_{T_n}(t) \|_{H^1}$ are uniformly in $t$ controlled by the mass and energy of $w_{T_n}$, and hence uniformly in $n$ bounded (see the more detailed argument in the proof of Assertion 2), say by $M$.  Thus, $\|u(t)\|_{L_t^\infty H_x^1} \leq M$.  The next goal is to confirm that $u(t)$ solves NLS in the sense of distributions (Assertion 5 below), and thus we can appeal to the Lemma \ref{L:NLS_uniqueness} giving that $u(t)$ is a strong solution to NLS and satisfies the mass and energy conservation laws.  In order to do this, we need Assertion 4 below.

\medskip

\noindent\textbf{Assertion 4}.  For the same sequence $T_n$ in Assertion 3, we have that for all $t$, $w_{T_n}(t)\to u(t)$ strongly in $L_{\text{loc}}^{p+1}$.

\medskip

\noindent\textit{Proof of Assertion 4}.  
Fix $t\in \mathbb{R}$, and pick any subsequence $T_{n_k}$ of $T_n$.  By the Rellich theorem, there exists a subsequence $T_{\tilde n_k}$ of $T_{n_k}$ such that $w_{T_{\tilde n_k}}(t)$ converges strongly in $L^{p+1}_{\text{loc}}$.  Since $w_{T_{\tilde n_k}}(t)$ converges weakly to $u(t)$ in $H^1$, we must have that $w_{T_{\tilde n_k}}(t) \to u(t)$ strongly in $L^{p+1}_{\text{loc}}$.  It follows\footnote{Recall the basic fact:  Given $h$ and a sequence $h_n$ in a metric space, suppose that for every subsequence $h_{n_k}$ of $h_n$, there exists a subsequence $h_{\tilde n_k}$ of $h_{n_k}$ such that $h_{\tilde n_k} \to h$.  Then $h_n \to h$, i.e. the original sequence converges to $h$.} that  $w_{T_n}(t) \to u(t)$ in $L^{p+1}_{\text{loc}}$.

\medskip

Now that we know that we have one distinguished sequence $T_n$ such that $\forall \; t$, $w_{T_n}(t) \to u(t)$ weakly in $H^1$ and strongly in $L^{p+1}_{\text{loc}}$, we can prove:

\medskip

\noindent \textbf{Assertion 5}.  $u$ solves NLS in the distributional sense.

\smallskip

\noindent \textit{Proof of Assertion 5}.
 Let $\phi\in C_c^\infty(\mathbb{R}^{n+1})$ be a test function.  We must show that (with $\la \cdot, \cdot \ra$ denoting the $H_x^1$--$H_x^{-1}$ pairing)
\begin{equation}
\label{E:u_solves_NLS}
\int_t \la u(t), -i\partial_t\phi(t) + \Delta\phi(t) \ra \, dt+ \int_t \la |u(t)|^{p-1}u(t),\phi(t) \ra \, dt =0 \, .
\end{equation}
We know that for each $n$,
\begin{equation}
\label{E:wn_solves_NLS}
\int_t \la w_{T_n}(t), -i\partial_t\phi(t) + \Delta\phi(t) \ra \, dt+ \int_t \la |w_{T_n}(t)|^{p-1}w_{T_n}(t),\phi(t) \ra \, dt =0 \, .
\end{equation}
For each fixed $t$, we have
 $$ \la w_{T_n}(t), -i\partial_t\phi(t) + \Delta\phi(t) \ra  \to  \la u(t), -i\partial_t\phi(t) + \Delta\phi(t) \ra \text{ as }n\to \infty$$
since $w_{T_n}(t) \to u(t)$ weakly; also we have
$$ \la |w_{T_n}(t)|^{p-1}w_{T_n}(t),\phi(t) \ra  \to  \la |u(t)|^{p-1}u(t),\phi(t) \ra \text{ as }n\to \infty $$
since $w_{T_n}(t) \to u(t)$ strongly in $L^{p+1}_{\text{loc}}$.  By dominated convergence (using that $\|w_{T_n}\|_{L_t^\infty H_x^1} \leq c$ independent of $n$), we can send $n\to \infty$ in \eqref{E:wn_solves_NLS} to obtain \eqref{E:u_solves_NLS}.  \textit{End proof of Assertion 5}.

\medskip

At this point we know that $u\in L_t^\infty H_x^1$ solves NLS as a distribution and thus by Lemma \ref{L:NLS_uniqueness}, $u\in C(\mathbb{R}; H_x^1)$ and $u$ satisfies the mass and energy conservation laws. 

It follows easily from Assertion 3 that for each $t\in (-\infty,+\infty)$, $U(-t)w_{T_n}(t)$ converges weakly to $U(-t)u(t)$ in $H^1$ as $n\to +\infty$.  We also have:

\medskip

\noindent\textbf{Assertion 6}.  $U(-t)u(t) \to v(0)$ weakly in $H^1$ as $t\to +\infty$.

\medskip

\noindent\textit{Proof of Assertion 6}. As in Assertion 1 we write
\begin{equation}
\label{E:A}
\langle U(-t)w_{T_{n}}(t), \psi_{j}\rangle=\langle v(0), \psi_{j}\rangle-i\int_{T_{n}}^{t}\langle U(-\sigma)\{|w_{T_{n}}|^{p-1}w_{T_{n}}(\sigma)\}, \psi_{j} \rangle d\sigma.
\end{equation}
By the above Arzela-Ascoli argument, we know that $\la U(-t)w_{T_n}(t),\psi_j \ra \to \la U(-t)u(t),\psi_j \ra$ \emph{uniformly in $t$} as $n\to +\infty$.  Let $\epsilon>0$.  It follows that $\exists \; N$ such that $n\geq N$ implies that 
\begin{equation}
\label{E:B}
\forall \; t, \qquad |\la U(-t)w_{T_n}(t),\psi_j \ra - \la U(-t)u(t),\psi_j \ra| < \frac{\epsilon}{2}
\end{equation}
But by Lemma \ref{L:NLSwkdecay}, for any $n$, 
\begin{equation}
\label{E:C}
\forall \; t\geq T_n, \quad \left| \int_{T_n}^t \la U(-\sigma)|w_{T_n}(\sigma)|^{p-1}w_{T_n}(\sigma), \psi_j \ra \, d\sigma \right| \lesssim T_n^{-\beta}
\end{equation}
Fix $n\geq N$ such that $T_n^{-\beta} \leq \frac{\epsilon}{2}$.  By combining \eqref{E:A}, \eqref{E:B}, and \eqref{E:C}, we obtain that for all $t\geq T_n$,
$$| \la U(-t)u(t),\psi_j\ra -\la v(0),\psi_j \ra | \leq \epsilon$$
\textit{End proof of Assertion 6}.

\medskip
Now we are ready to finish the proof.  The key facts that we have established thus far are:  (1) $u(t)\in C(\mathbb{R}; H_x^1)$ solves NLS and satisfies mass and energy conservation; (2) for fixed $t$, $w_{T_n}(t) \to u(t)$ weakly in $H_x^1$ as $n\to \infty$; (3) $U(-t)u(t) \to v(0)$ weakly in $H_x^1$ as $n\to \infty$.  The main remaining goal is to upgrade the weak convergence in (3) to strong convergence.  In the process, we will end up upgrading the weak convergence in (2) to strong convergence.

Since $w_{T_n}(t)$ converges weakly to $u(t)$ in $H^1$ as $n\to \infty$, and using the $L^2$ conservation of NLS and the free Schr\"odinger equation, we have
$$\|u(t)\|_{L^2} \leq \liminf_{n\to +\infty} \|w_{T_n}(t) \|_{L^2} \leq \limsup_{n\to +\infty} \|w_{T_n}(t) \|_{L^2} = \|v(0)\|_{L^2}$$
However, since $U(-t)u(t)$ converges weakly to $v(0)$ in $H^1$ as $t\to +\infty$, we have
$$\|v(0)\|_{L^2} \leq \liminf_{t\to \infty} \|U(-t)u(t)\|_{L^2} \leq \limsup_{t\to \infty} \|U(-t)u(t)\|_{L^2} = \|u(t)\|_{L^2}$$
Combining the two inequalities, we learn that all of the inequalities are equalities and thus that $U(-t)u(t) \to v(0)$ strongly in $L^2$ as $t\to +\infty$ and that for each $t\in (-\infty,+\infty)$, $w_{T_n}(t)\to u(t)$ strongly in $L^2$. To prove the above statements we use the fact that if $f_{j} \rightarrow f$ weakly in $L^2$ and $\|f_{j}\|_{L^2} \rightarrow \|f\|_{L^2}$ then $f_{j} \rightarrow f$ strongly in $L^2$. 
By interpolation inequalities, the fact that these limits are strong in $L^2$ and bounded in $H^1$, we learn that these limits are strong in $L^{p+1}$.

Since $U(-t)u(t)$ converges to $v(0)$ weakly in $H^1$ as $t\to +\infty$, we have
\begin{equation}
\label{E:wkbd101}
\|\nabla v(0)\|_{L^2}^2 \leq \liminf_{t\to +\infty} \|\nabla \; U(-t)u(t)\|_{L^2}^2 = \liminf_{t\to +\infty} \|\nabla u(t)\|_{L^2}^2
\end{equation}
Since for each $t\in (-\infty,+\infty)$, $w_{T_n}(t)$ converges weakly to $u(t)$ in $H^1$, we have
\begin{align*}
E[u(t)] &= \frac12 \|\nabla u(t)\|_{L^2}^2 + \frac1{p+1}\|u(t)\|_{L^{p+1}}^{p+1}\\
&\leq \liminf_{n\to +\infty} \; \frac12\|\nabla w_{T_n}(t)\|_{L^2}^2 \; + \frac1{p+1}\|u(t)\|_{L^{p+1}}^{p+1}\\
&= \liminf_{n\to +\infty} \left( E[w_{T_n}(t)] - \frac1{p+1}\|w_{T_n}(t)\|_{L^{p+1}}^{p+1}\right)  + \frac1{p+1}\|u(t)\|_{L^{p+1}}^{p+1}\\
&= \liminf_{n\to +\infty}  E[w_{T_n}(T_n)] \,,
\intertext{where we have used that $w_{T_n}(t)$ converges strongly to $u(t)$ in $L^{p+1}$ and the conservation of energy for $w_{T_n}$.  Since $w_{T_n}(T_n)=v(T_n)$, we have,}
&=  \liminf_{n\to +\infty} \left( \frac12\|\nabla v(T_n)\|_{L^2}^2 - \frac1{p+1}\|v(T_n)\|_{L^{p+1}}^{p+1} \right)\\
&= \frac12\|\nabla v(0)\|_{L^2}^2 \,,
\end{align*}
where we have used that the free Schr\"odinger evolution is unitary on $\dot H^1$ and the space-time decay estimate.  Thus, we have
\begin{equation}
\label{E:wkbd100}
E[u(t)]\leq \frac12\|\nabla v(0)\|_{L^2}^2
\end{equation}
Since $U(-t)u(t)$ converges strongly to $v(0)$ in $L^2$ as $t\to +\infty$, we have that $\|u(t)-v(t)\|_{L^2} \to 0$ as $t\to +\infty$.  By interpolation and the boundedness of $u(t)$ and $v(t)$ in $H^1$, we have that $\|u(t)-v(t)\|_{L^{p+1}} \to 0$ as $t\to +\infty$.  But the space-time decay estimate for free Schr\"odinger propagation implies that $v(t)\to 0$ in $L^{p+1}$ and thus $u(t)\to 0$ in $L^{p+1}$ as $t\to +\infty$.  By \eqref{E:wkbd100}, we have
$$\limsup_{t\to +\infty}  \|\nabla u(t)\|_{L^2}^2 \leq \|\nabla v(0)\|_{L^2}^2 \,.$$
Combining this with \eqref{E:wkbd101}, we see that the inequalities must be equalities and thus $U(-t)u(t) \to v(0)$ strongly in $H^1$.

\end{document}